\documentclass[10pt]{amsart}
\usepackage{amsfonts}
\usepackage{amssymb}
\usepackage{amsthm}
\usepackage{amsmath}
\usepackage{graphicx}
\theoremstyle{plain}
\newtheorem{theorem}{Theorem}[section]
\theoremstyle{remark}
\newtheorem*{remark}{Remark}
\usepackage[all]{xy}
\author{C. Robinson Tompkins}
\title{Latin Square Thue-Morse Sequences are Overlap-Free}
\date{}
\begin{document}

\begin{abstract}
   We define a morphism based upon a Latin square that generalizes the
   Thue-Morse morphism. We prove that fixed points of this
   morphism are overlap-free sequences, generalizing results
   of Allouche - Shallit and Frid.
\end{abstract}

\maketitle

\section{Introduction}

In his 1912 paper, Axel Thue introduced the first binary sequence
that does not contain an overlap \cite{th}. It is now called the
Thue-Morse sequence:
\begin{displaymath}
  0110 1001 1001 0110 1001 0110 0110 1001\ldots.
\end{displaymath}
An overlap is a string of letters in the form
$c\mathbf{x}c\mathbf{x}c$ where $c$ is a single letter and
$\mathbf{x}$ is finite string that is potentially empty. Overlaps
begin with a square, namely $\mathbf{w}\mathbf{w}$ where
$\mathbf{w}=c\mathbf{x}$ as given above. It is easy to observe, as
Thue did, that any binary string of four or more letters must contain
a square.

There are several ways to define the Thue-Morse sequence \cite{al2}.
We will derive it as a fixed point of a morphism. Let $\Sigma$ be an
alphabet and let $\Sigma^*\cup\Sigma^\omega$ be the set of all finite or infinite strings over $\Sigma$. A morphism is a mapping
\begin{displaymath}
  h:\Sigma^*\cup\Sigma^\omega\rightarrow\Sigma^*\cup\Sigma^\omega
\end{displaymath}
that obeys the identity $h(xy)=h(x)h(y)$, for $x$ a finite string
and $y\in\Sigma^*\cup\Sigma^\omega$ \cite[p. 8]{al}.

By \cite[p. 16]{al}, define the Thue-Morse morphism on
$\Sigma=\{0,1\}$ as
\begin{equation}
  \mu(t)=\left\{
  \begin{array}{ll}
    01,&\textrm{for}\;t=0\\
    10,&\textrm{for}\;t=1
  \end{array}\right..
\end{equation}
The sequence found by applying 0 to the the $n$th iterate 
of $\mu$  converges to the Thue-Morse sequence, denoted 
$\mu^\omega(0)$, which of course is infinite. In particular,
\begin{displaymath}
  \begin{array}{rcl}
  \mu(0)&=&01\\
  \mu^2(0)&=&\mu(\mu(0))=\mu(01)=\mu(0)\mu(1)=0110\\
  \mu^3(0)&=&\mu(\mu^2(0))=\mu(0110)=01101001\\
  &\vdots&\\
  \mu^\omega(0)&=&0110 1001 1001 0110 1001 0110 0110 1001\ldots.
  \end{array}
\end{displaymath}
Notice that $\mu^\omega(\mu(0))=\mu^\omega(0)$ and
$\mu(\mu^\omega(0))=\mu^\omega(0)$. This second observation says
that the Thue-Morse sequence is a fixed point of $\mu$ \cite[p.
10]{al}.

We can identify the binary alphabet of the Thue-Morse sequence with
$\mathbb{Z}/2\mathbb{Z}$ the integers modulo 2. It is natural to then
generalize it to $\mathbb{Z}/n\mathbb{Z}$, by considering the alphabet
$\Sigma=\{0,1,\ldots,n-1\}$, and for $i\in\Sigma$, defining the
morphism
\begin{displaymath}
  \phi_n(i)=\overline{i+0}\;\overline{i+1}\ldots\overline{i+(n-1)},
\end{displaymath}
where $\overline{i}$ is the residue modulo $n$. Notice that for
$\Sigma=\{0,1\}$, $\phi_2(i)=\mu(i)$. In 2000, Allouche and Shallit
proved that $\phi_n^\omega$ is overlap-free \cite{al3}.

In this paper, we generalize $\phi_n$, which is based on the Cayley
table of $\mathbb{Z}/n\mathbb{Z}$, to Latin squares of arbitrary finite size
$n$. We define our morphism based the Latin square, and prove that
the fixed point of the Latin square morphism is an overlap-free
sequence. Note that the Cayley table for $\mathbb{Z}/n\mathbb{Z}$ is a Latin
square, but not every Latin square is a Cayley table.

\section{Latin Square Morphisms produce Tilings}

Allouche and Shallit's morphism can be seen as a mapping of $i$ to
the $i^\mathrm{th}$ row (that begins with $i$) of the Cayley table
for $\mathbb{Z}/n\mathbb{Z}$. For example when $n=3$, we have

\begin{displaymath}
  \begin{array}{ccccc}
    &\phi_3&&&\\
    0&\rightarrow&0&1&2\\
    1&\rightarrow&1&2&0\\
    2&\rightarrow&2&0&1
  \end{array}
\end{displaymath}
This suggests a natural generalization to any Latin square.

Begin with a generic alphabet of $n$ letters, which we may assume
to be $\{1,2,\ldots,n\}$. Recall that a Latin square $\mathcal{L}$
is an $n\times n$ table with $n$ different letters such that each
letter occurs only once in each column and only once in each row. We
will concern ourself with the Latin squares in which the first
column retains the natural order of our alphabet $(1,2,\ldots,n)$.
For $n=3$, there are two such Latin squares. The one that does not
come from $\mathbb{Z}/3\mathbb{Z}$ directly is
\begin{displaymath}
  \left[
  \begin{array}{ccc}
    1&3&2\\
    2&1&3\\
    3&2&1
  \end{array}
  \right].
\end{displaymath}

Let $\mathcal{L}_t$ denote the $t^\mathrm{th}$ row of our Latin
square $\mathcal{L}$. For each $t\in\Sigma$ we define the Latin
square morphism by $\ell(t)=\mathcal{L}_t$. For example we can use
the above Latin square for $n=3$ to define the following morphism,
\begin{displaymath}
  \ell(t)=\left\{
  \begin{array}{ll}
    132,&\textrm{for}\;t=1\\
    213,&\textrm{for}\;t=2\\
    321,&\textrm{for}\;t=3
  \end{array}\right.
\end{displaymath}
Given any $t\in\Sigma$, $\ell(t),\ell^2(t),\ell^3(t),\ldots$
converges to a sequence $\ell^\omega(t)$, which is a fixed point of
the morphism $\ell$. So,
\begin{equation}
  \ell(\ell^\omega(t))=\ell^\omega(t)
\end{equation}
In fact every fixed point of $\ell$ is of the form $\ell^\omega(t)$
for some $t\in\Sigma$ \cite[p. 10]{al}.

Express the sequence as $\ell^\omega(t_1) = t_1t_2t_3\ldots$, so
\begin{displaymath}
    \ell^\omega(t_1)=\ell(\ell^\omega(t_1))=\ell(t_1t_2t_3\ldots)=\ell(t_1)\ell(t_2)\ell(t_3)\ldots
\end{displaymath}
Thus, we have a tiling of our sequence (and of the natural numbers)
by the rows of our Latin square $\mathcal{L}$. Again, in terms of
our example where $n=3$ we have three tiles 132, 213, and 321 and so
\begin{displaymath}
  \ell^\omega(1)=132321213321213132\ldots=|132|321|213|321|213|132|\ldots.
\end{displaymath}
Now, consider the subsequence created by taking the first letter of
each tile. Notice that this sequence is in fact our original
sequence. Thus our sequence contains itself as a subsequence. These
two observations, our sequence as a tiling and our sequence equaling
a subsequence of itself, will be critical for the proof of our main
result.

\section{Overlap-Free Latin Square Sequences}

In this section we prove our main result.

\begin{theorem}
  Let $\Sigma=\{1,2,\ldots,n\}$, and let $\mathcal{L}$ be an
  $n\times n$ Latin square using the letters from $\Sigma$,
  with the first column in its natural order. For an arbitrary
  $t\in\Sigma$, let $\mathcal{L}_t$ denote the row of
  $\mathcal{L}$ corresponding to $t$ in the first column. If we
  define the Latin square morphism as
  \begin{displaymath}
    \ell(t)=\mathcal{L}_t,
  \end{displaymath}
  then we have that for any $t\in\Sigma$, $\ell^\omega(t)$ is an
  overlap-free sequence.
\end{theorem}

\begin{remark}
The Latin square for $n=3$ above can be seen to be the Cayley table
for $\mathbb{Z}/3\mathbb{Z}$ with the last two columns transposed. Frid has
shown that all morphisms based upon such Latin squares for
$\mathbb{Z}/n\mathbb{Z}$ produce overlap-free sequences as their fixed points
\cite{fr}. Of course not every Latin square comes from a group
Cayley table. For an example of a Latin square that is not a group
Cayley table see below \cite[p. 27]{de}.
\begin{displaymath}
\left[\begin{array}{cccccc}
1&2&3&4&5&6\\
2&1&6&3&4&5\\
3&4&5&2&6&1\\
4&5&1&6&2&3\\
5&6&4&1&3&2\\
6&3&2&5&1&4
\end{array}\right]
\end{displaymath}
\end{remark}

\begin{proof}
  Let $\ell^\omega(t_1) = t_1t_2t_3\ldots$ so the $j^\mathrm{th}$ letter in the sequence is $t_j$.
  Similarly, the $m^\mathrm{th}$ tile in the
  sequence is $T_m$. We will be also using the notion of length of a
  string of letters, meaning the number of letters in a string. For an
  arbitrary string $w$ the length of $w$ will be denoted $|w|$.
  Use $r$ to denote the location of $t_j$ on its tile
  $T_m$, so $j=(m-1)n+r$ with $|T_m|=n$ and $r\in\{1,2,\ldots,n\}$.

  Assume for a contradiction that $\ell^\omega(t_1)$ contains an
  overlap; moreover that $c\mathbf{x}c\mathbf{x}c$ is the shortest
  overlap in $\ell^\omega(t_1)$. Write
  $\ell^\omega(t_1)=Ac\mathbf{x}c\mathbf{x}cB$, where $c$ is a single
  letter, $\mathbf{x}$ is a finite string with $|c\mathbf{x}|\geq
  n$, $A$ is a finite string, and $B$ is the infinite tail of our
  sequence. We have that $|c\mathbf{x}|\geq n$ (bound by the length
  of the tiles) because each tile is a permutation of
  $1,2,\ldots,n$, and we cannot have two of the three copies of $c$
  contained in one tile. Our subscripts place this overlap
  in our sequence. For $i\in\{1,2,3\}$, let $j_i$ denote the
  subscript of the $i^\mathrm{th}$ $c$. Thus,
  \begin{equation}
    \begin{array}{rcl}
      A&=&t_1\cdots t_{j_1-1}\\
      c&=&t_{j_1}=t_{j_2}=t_{j_3}\\
      \mathbf{x}&=&t_{j_1+1}\cdots t_{j_2-1}=t_{j_2+1}\cdots t_{j_3-1}\\
      B&=&t_{j_3+1}t_{j_3+2}t_{j_3+3}\cdots,
    \end{array}
  \end{equation}

  Our argument proceeds as follows: there are two cases
  $|c\mathbf{x}|\not\equiv0\pmod{n}$ and
  $|c\mathbf{x}|\equiv0\pmod{n}$. In the first case we use the
  fact that we have a tiling of $\ell^\omega(t_1)$ by the rows of a
  Latin square, to show that the overlap $c\mathbf{x}c\mathbf{x}c$ is
  not possible. In the
  second case, when $|c\mathbf{x}|\equiv0\pmod{n}$, we argue based
  upon the fact that $\ell^\omega(t_1)$ contains itself as a subsequence
  that the existence of the overlap $c\mathbf{x}c\mathbf{x}c$
  leads to the existence of a shorter overlap, and thus a
  contradiction.

  \subsection{Case 1: $|c\mathbf{x}|\not\equiv 0
  \pmod{n}$}

  For each $i\in\{1,2,3\}$, let $r_i\in\{1,2,\ldots,n\}$ such that
  $r_i\equiv j_i\pmod{n}$. In other words $t_{j_i}$ is the
  $r_i^{\phantom{i}\mathrm{th}}$ letter in its tile in
  $\ell^\omega(t_1)$. Also, we will refer to the tile containing $t_{j_i}$ as
  $T_{m_i}$. It is now possible to write the length of
  $c\mathbf{x}$ as $|c\mathbf{x}|\equiv r_2-r_1\equiv
  r_3-r_2\pmod{n}$. So,
  \begin{equation}
    r_3\equiv 2r_2 - r_1\pmod{n}.
  \end{equation}

  \subsubsection{Six Cases}

  Since $r_2-r_1\equiv |c\mathbf{x}|\not\equiv0\pmod{n}$
  there are two main cases that we will first consider: $r_1<r_2$
  and $r_2<r_1$. However, for the explicit details of our
  conclusions we will consider all six of the following possibilities depending on
  the value of $r_3$,
  \begin{displaymath}
    \begin{array}{c}
      r_3=2r_2-r_1\longleftrightarrow\left\{
      \begin{array}{c}
        r_1<r_2<r_3\\
        r_3<r_2<r_1
      \end{array}\right.\\
      r_3=2r_2-r_1-n\longleftrightarrow\left\{
      \begin{array}{c}
        r_1\leq r_3<r_2\\
        r_3<r_1<r_2
      \end{array}\right.\\
      r_3=2r_2-r_1+n\longleftrightarrow\left\{
      \begin{array}{c}
        r_2<r_1\leq r_3\\
        r_2<r_3<r_1
      \end{array}\right.
    \end{array}
  \end{displaymath}
  The equalities on the left arise out of equation (4) and the fact
  that the integer $2r_2-r_1$ satisfies, $-n\leq 2r_2-r_1\leq 2n$.
  This means that $r_3$ is the element in the set
  $\{2r_2-r_1+n,2r_2-r_1,2r_2-r_1-n\}$ that lies in the interval
  $0<r_3\leq n$. Notice that $r_3=2r_2-r_1$ in both cases when
  $r_1<r_2$ and $r_2<r_1$.

  \subsubsection{$G$ and the beginning of each $c\mathbf{x}$}

  When $r_1<r_2$, we pick $G\subset\Sigma$ to be the last $r_2-r_1$
  letters in $T_{m_1}$ such that $G$ has no specific order and
  $G\neq\emptyset$. Of course, the remainder of the letters in
  $T_{m_1}$ are in $\overline{G}$, the complement of $G$. Notice
  that this puts $c=t_{j_1}\in\overline{G}$. By equating the letters
  in $T_{m_1}$ with the corresponding letters in
  $t_{j_2}\mathbf{x}t_{j_3}$, we find that the last $n-r_2+1$ letters
  of $T_{m_2}$
  (starting with $c=t_{j_2}$) are in $\overline{G}$. Also, we find that
  the first $r_2-r_1$ letters of $T_{m_2+1}$ are $G$.

  When $r_2<r_1$, we pick $G\subset\Sigma$ to be the last $r_1-r_2$
  letters in $T_{m_2}$ such that $G$ has no specific order and
  $G\neq\emptyset$. Obviously, the remainder of letters in $T_{m_2}$
  must be those that make up $\overline{G}$ again placing
  $c=t_{j_2}\in\overline{G}$. By equating the letters in $T_{m_2}$
  with the corresponding letters in $t_{j_1}\mathbf{x}t_{j_2}$ we find
  that the last $n-r_1+1$ letters of $T_{m_1}$ (starting with $c=t_{j_1}$)
  are in $\overline{G}$. Also, we find that the first $r_1-r_2$
  letters of $T_{m_1+1}$ are $G$.

  We have discussed the appearance of $G$ and its complement
  $\overline{G}$ in the beginning of each $c\mathbf{x}$. So, we set
  forth to describe $G$ and $\overline{G}$ at the end of each
  $c\mathbf{x}$.

  \subsubsection{Following $G$ through the overlap}

  It is a basic observation that because each tile is a permutation
  of the letters in $\Sigma$, each tile can be partitioned into $G$
  and its complement $\overline{G}$. It is fundamental to our
  argument that because of the equality $t_{j_1}\mathbf{x}t_{j_2}=
  c\mathbf{x}c=t_{j_2}\mathbf{x}t_{j_3}$, the letters in $G$ form a
  contiguous collection of elements in each tile involved in our
  overlap excluding $T_{m_i}$ (each of which will need further description),
  either the beginning or the ending of each tile. The idea
  involved in following $G$ through the overlap is quite simple,
  we illustrate it in one particular case $r_1<r_2<r_3$.

  We have explicitly described the location of $G$ at
  the beginning of each $c\mathbf{x}$. We will now use our example
  $r_1<r_2<r_3$ to show to the reader how the tiling of our sequence
  can be used to find the location of $G$ at the end of each
  $c\mathbf{x}$. In doing so, we will refer to Figure 1.

  In Figure 1, we have displaced the overlap from our sequence (represented by the
  continuous solid horizontal line). We have also split
  our overlap in half leaving $T_{m_2}$ intact for equality
  purposes. We have placed $t_{j_1}\mathbf{x}t_{j_2}$ over
  $t_{j_2}\mathbf{x}t_{j_3}$ with $t_{j_1}$ directly over $t_{j_2}$
  and $t_{j_2}$ directly over $t_{j_3}$ so that we can see equality
  of terms simply by looking straight up or straight down (displayed
  by vertical arrows). The set of letters $G$ is represented by a
  horizontal solid line above and below our sequence line, and the set of letters
  $\overline{G}$ is represented by horizontal dotted lines above
  and below the sequence line. Also,
  notice that we have drawn in the edges of the tiles with smaller
  vertical black lines.
  \begin{displaymath}
    \begin{xy}
        (-50,5);(60,5) **@{-},
        (-50,4);(-50,6) **@{-},
        (-40,4);(-40,6) **@{-},
        (-30,4);(-30,6) **@{-},
        (-20,4);(-20,6) **@{-},
        (-10,4);(-10,6) **@{-},
        (0,4);(0,6) **@{-},
        (10,4);(10,6) **@{-},
        (20,4);(20,6) **@{-},
        (30,4);(30,6) **@{-},
        (40,4);(40,6) **@{-},
        (50,4);(50,6) **@{-},
        (60,4);(60,6) **@{-},
        (-48,3);(-48,7) **@{-},
        <-47mm,9mm>*\txt\small{$t_{j_1}$},
        (54,3);(54,7) **@{-},
        <55mm,9mm>*\txt\small{$t_{j_2}$},
        (-50,5);(-42,5) **\dir2{.},
        (-42,5);(-40,5) **\dir2{-},
        (-40,5);(-32,5) **\dir2{.},
        (-32,5);(-30,5) **\dir2{-},
        (-30,5);(-22,5) **\dir2{.},
        (-22,5);(-20,5) **\dir2{-},
        (38,5);(40,5) **\dir2{-},
        (40,5);(48,5) **\dir2{.},
        (48,5);(50,5) **\dir2{-},
        (54,5);(60,5) **\dir2{.},
        <-41.5mm,7.5mm>*\txt\small{$G$},
        <-35.5mm,7.82mm>*\txt\small{$\overline{G}$},
        <48.5mm,7.5mm>*\txt\small{$G$},
        <57mm,2mm>*\txt\small{$\overline{G}$},
        (-52,-5);(58,-5) **@{-},
        (-52,-6);(-52,-4) **@{-},
        (-42,-6);(-42,-4) **@{-},
        (-32,-6);(-32,-4) **@{-},
        (-22,-6);(-22,-4) **@{-},
        (-12,-6);(-12,-4) **@{-},
        (-2,-6);(-2,-4) **@{-},
        (8,-6);(8,-4) **@{-},
        (18,-6);(18,-4) **@{-},
        (28,-6);(28,-4) **@{-},
        (38,-6);(38,-4) **@{-},
        (48,-6);(48,-4) **@{-},
        (58,-6);(58,-4) **@{-},
        (-48,-7);(-48,-3) **@{-},
        <-48mm,-9mm>*\txt\small{$t_{j_2}$},
        (54,-7);(54,-3) **@{-},
        <54mm,-9mm>*\txt\small{$t_{j_3}$},
        (-48,-5);(-42,-5) **\dir2{.},
        (-42,-5);(-40,-5) **\dir2{-},
        (-40,-5);(-32,-5) **\dir2{.},
        (-32,-5);(-30,-5) **\dir2{-},
        (-30,-5);(-22,-5) **\dir2{.},
        (38,-5);(40,-5) **\dir2{-},
        (40,-5);(48,-5) **\dir2{.},
        (48,-5);(50,-5) **\dir2{-},
        (50,-5);(58,-5) **\dir2{.},
        (-41,-4);(-41,4) **@{-}
        ?<*@{<} ?>*@{>},
        (-36,-4);(-36,4) **@{-}
        ?<*@{<} ?>*@{>},
        (-31,-4);(-31,4) **@{-}
        ?<*@{<} ?>*@{>},
        (-26,-4);(-26,4) **@{-}
        ?<*@{<} ?>*@{>},
        (39,-4);(39,4) **@{-}
        ?<*@{<} ?>*@{>},
        (44,-4);(44,4) **@{-}
        ?<*@{<} ?>*@{>},
        (49,-4);(49,4) **@{-}
        ?<*@{<} ?>*@{>},
        <15mm,0mm>*\txt\small{$\bullet$},
        <20mm,0mm>*\txt\small{$\bullet$},
        <25mm,0mm>*\txt\small{$\bullet$},
    \end{xy}
  \end{displaymath}
  \begin{center}
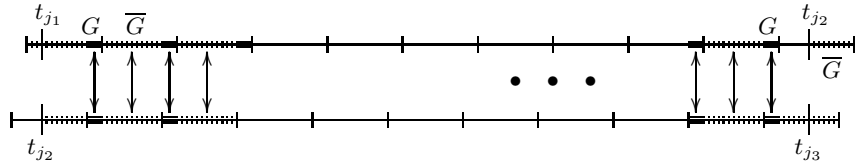
     Figure 1: The situation when $r_1<r_2<r_3$.
  \end{center}

  Now notice that by using the tiles we can equate letters in
  $t_{j_1}\mathbf{x}t_{j_2}$ with $t_{j_2}\mathbf{x}t_{j_3}$ all the
  way through the overlap. Since we know that $G$ occurs in the
  first $r_2-r_1$ letters of $T_{m_2+1}$, then $\overline{G}$ is
  the last $n-(r_2-r_1)$ letters of $T_{m_2}+1$. This causes
  $\overline{G}$ to be the first $n-(r_2-r_1)$ letters of
  $T_{m_1+1}$, and thus $G$ appears in the last
  $r_2-r_1$ letters of $T_{m_1+1}$. Thus we can conclude that $G$
  occurs in the last $r_2-r_1$ letters of all the tiles in
  $t_{j_1}\mathbf{x}t_{j_2}$ except for $T_{m_2}$. We can also
  conclude that $G$ occurs in the first $r_2-r_1$ letters of all the
  tiles in $t_{j_2}\mathbf{x}t_{j_3}$ up through $T_{m_3-1}$. We can approach every case by
  the same process.

  \subsubsection{$G$ and how each $c\mathbf{x}$ ends}

  We now will explain the conclusions for the six possible
  cases that we defined earlier, leaving the actual drawing to the reader.

  Case $r_1<r_2<r_3$ (as seen in Figure 1). After
  we follow $G$ through the overlap, we find that $G$ occurs in the
  first $r_2-r_1$ letters of $T_{m_3}$. Recall $r_3=2r_2-r_1$. So, we have that the next
  $r_3-(r_2-r_1)=r_2$ letters of $T_{m_3}$ are not
  in $G$. Notice that the size of $G$, $r_2-r_1$, added to $r_2$ make
  up all of $r_3$. This places the boundary between $T_{m_2-1}$ and $T_{m_2}$
  exactly in line with the end of $G$ in $T_{m_3}$ and the beginning
  of $\overline{G}$. We then equate the first letters in $T_{m_3}$ with those in
  $T_{m_2}$ to find that $G$ occurs nowhere in $T_{m_2}$. So now, we
  have described $T_{m_2}$ fully. Earlier we defined $G$ such that
  $\overline{G}$ occurred from $t_{j_2}$ to the end of the tile, and
  we have just shown that the first $r_2$ letters of $T_{m_2}$ (which
  includes $t_{j_2}$) must be in $\overline{G}$. So $G$
  does not appear in anywhere in $T_{m_2}$, and since
  $G\neq\emptyset$, we must have a contradiction.

  Cases $r_1\leq r_3<r_2$ and $r_3<r_1<r_2$.
  After we follow $G$ through the overlap, we find that $G$ occurs
  in the first $r_2-r_1$ letters of $T_{m_3-1}$. So, $\overline{G}$
  occurs in the final $n-(r_2-r_1)$ letters of $T_{m_3-1}$ causing the
  first $n-(r_2-r_1)$ letters of $T_{m_2}$ to be $\overline{G}$.
  Notice that $r_2=[n-(r_2-r_1)]+r_3$. So the boundary between
  $\overline{G}$ and $G$ in $T_{m_2}$ coincides with the boundary
  between $T_{m_3-1}$ and $T_{m_3}$. This means that $t_{j_2}\in G$,
  but we assumed that $c\notin G$ earlier which is a contradiction.

  Case $r_3<r_2<r_1$. After we follow $G$ through the
  overlap, we find that $G$ occurs in the last $r_1-r_2$ letters of
  $T_{m_3-1}$. This causes $G$ to occur in the first $r_1-r_2$ letters
  of $T_{m_2}$ by equality of $t_{j_1}\mathbf{x}t_{j_2}$ and
  $t_{j_2}\mathbf{x}t_{j_3}$. To describe the remaining letters of
  $T_{m_2}$ up to and including $t_{j_2}$ consider
  $r_2-(r_1-r_2)=r_3$. So $\overline{G}$ occurs in the next $r_3$
  letters after $G$. Thus we have
  that $G$ is repeated twice in $T_{m_2}$ so we have our
  contradiction.

  Cases $r_2<r_1\leq r_3$ and $r_2<r_3<r_1$. After we
  follow $G$ through the overlap we find that $G$ occurs in the
  first $r_1-r_2$ letters of $T_{m_2-1}$. This causes $\overline{G}$
  to occur in the final $n-(r_1-r_2)$ letters of $T_{m_2-1}$ and thus
  the first $n-(r_1-r_2)$ letters of $T_{m_3}$. Since
  $r_2=r_3-[n-(r_1-r_2)]$, we see that the left boundary of
  $T_{m_2}$ coincides with the right boundary of these first
  $n-(r_1-r_2)$ letters of $T_{m_3}$. In particular, this means that
  the last $r_1-r_2$ letters of $T_{m_3}$, which include $c$, are in
  $G$. But, this contradicts the fact that $c\notin G$.

  \subsection{Case 2: $|c\mathbf{x}|\equiv 0 \pmod{n}$}

  We begin by considering some $\pi\in S_n$ the symmetric group on
  $n$ letters. Note that we may apply $\pi$ to any string by
  requiring $\pi$ to act on each individual letter, so
  $\pi(t_1t_2\ldots t_s)=\pi(t_1)\pi(t_2)\ldots\pi(t_s)$. Thus
  $\pi$ can be treated as a morphism. Moreover,
  $\pi:\Sigma^*\rightarrow\Sigma^*$ is an invertible map because $\pi\in
  S_n$. Thus $w\in\Sigma^*$ contains an overlap if and only if
  $\pi(w)\in\Sigma^*$ contains an overlap.

  Define the function
  $d_{(a,n)}:\mathbb{N}\rightarrow\mathbb{N}$ by $d_{(a,n)}(m)=(m-1)n+a$.
  Now if we let $M=(t_s)$ be a sequence, then define the sequence
  given by the function $D_{(a,n)}(M)$ to be the subsequence
  $(t_{d_{(a,n)}(s)})$ of $M$. So for $i\in\{1,2,\ldots,n\}$
  arbitrary we have that
  \begin{displaymath}
    D_{(i,n)}(\ell^\omega(t_1))=t_it_{i+n}t_{i+2n}\ldots.
  \end{displaymath}

  Define $\pi_i:\Sigma\rightarrow\Sigma$ with $\pi_i\in S_n$, such that if
  $\mathcal{L}_{t_1}=\{t_1,t_2,\ldots,t_i,\ldots,t_n\}$,
  $\pi_i(t_1)=t_i$. Recall that $\mathcal{L}_t$ refers to the
  $t^{\mathrm{th}}$ row of our Latin square $\mathcal{L}$.
  So we have that $\pi_i$ maps each letter in the
  first column of our Latin square, to the $i^\mathrm{th}$ letter of
  its corresponding row. Now, we want to show that
  $\pi_i(\ell^\omega(t))=D_{(i,n)}(\ell^\omega(t))$ for all $t\in\Sigma$. So take
  \begin{displaymath}
    \begin{array}{rcl}
      D_{(i,n)}(\ell^\omega(t_1))&=&D_{(i,n)}(\ell(\ell^\omega(t_1))\\
      &=&D_{(i,n)}(\ell(t_1)\ell(t_2)\ell(t_3)\cdots)\\
      &=&\pi_i(t_1)\pi_i(t_2)\pi_i(t_3)\cdots\\
      &=&\pi_i(\ell^\omega(t_1)).
    \end{array}
  \end{displaymath}
  Since $\pi_i\in S_n$ is invertible we can conclude that
  $D_{(i,n)}(\ell^\omega(t_1))$ contains an overlap if and only if
  $\ell^\omega(t_1)$ contains an overlap.

  Since $|c\mathbf{x}|\equiv 0\pmod{n}$
  pick $i\equiv j_1\equiv j_2\equiv j_3\pmod{n}$. By applying $D_{(i,n)}$ to (4) we
  obtain
  \begin{displaymath}
    D_{(i,n)}(\ell^\omega(t_1))=
    A_i\;t_{j_1}\;\mathbf{x}_i\;t_{j_2}\;\mathbf{x}_i\;t_{j_3}\;B_i
  \end{displaymath}
  where
  \begin{displaymath}
    \begin{array}{rcccl}
      A_i&=&D_{(i,n)}(A)&=&t_it_{i+n}t_{i+2n}\ldots,\\
      \mathbf{x}_i&=&D_{(i,n)}(\mathbf{x})&=&t_{j_1+n}t_{j_1+2n}\ldots
      t_{j_1+(m-1)n}\\
      &&&=&t_{j_2+n}t_{j_2+2n}\ldots
      t_{j_2+(m-1)n},\\
      B_i&=&D_{(i,n)}(B)&=&t_{j_3+n}t_{j_3+2n}t_{j_3+3n}\ldots,
    \end{array}
  \end{displaymath}
  and $m=|c\mathbf{x}|/n$.
  Observe that $D_{(i,n)}(\ell^\omega(t_1))$ contains a shorter
  overlap which implies that $\ell^\omega(t_1)$ also contains
  a shorter overlap, a contradiction of our assumption.
\end{proof}

\section{Acknowledgements}

I would like to thank Dr. Griff Elder, my research advisor for his
guidance. I would also like to thank Dr. Dan Farkas for introducing
me to the Thue-Morse sequence at Virginia Tech's Undergraduate
Research workshop in 2006 funded by the NSA and for the idea that
lead to the argument for $|c\mathbf{x}|\equiv 0\pmod{n}$. Lastly, I would
like to thank the referees and Professor Anca Muscholl for helping
 me correct my notation.

\end{document}